# Directional Differentiability of the Metric Projection in Uniformly Convex and Uniformly Smooth Banach Spaces

Jinlu Li

Department of Mathematics
Shawnee State University
Portsmouth, Ohio 45662
USA

**Abstract**

Let $X$ be a real uniformly convex and uniformly smooth Banach space and $C$ a nonempty closed and convex subset of $X$. Let $P_C: X \to C$ denote the (standard) metric projection operator. In this paper, we define the Gâteaux directional differentiability of $P_C$. We investigate some properties of the Gâteaux directional differentiability of $P_C$. In particular, if $C$ is a closed ball or a closed and convex cone (including proper closed subspaces), then, we give the exact representations of the directional derivatives of $P_C$. Finally, we define the concept of $p$-$q$ uniformly convex and uniformly smooth Banach spaces. We will prove that if $X$ is a $p$-$q$ uniformly convex and uniformly smooth Banach space, then for any nonempty closed and convex subset $C$ of $X$, $P_C$ is directionally differentiable on the whole space $X$. The results in this paper can be immediately applied to Hilbert spaces.

## 1. Introduction

Let $(X, \|\cdot\|)$ be a real uniformly convex and uniformly smooth Banach space with topological dual space $(X^*, \|\cdot\|_*)$. Let $C$ be a nonempty closed and convex subset of $X$. Let $P_C: X \to C$ denote the (standard) metric projection operator. For any $x \in X$, $P_C x \in C$ such that

$$\|x - P_C x\| \leq \|x - z\|, \text{ for all } z \in C. \tag{1.1}$$

$P_C x$ is called the metric projection of point $x$ onto $C$. $P_C x$ is considered as the best approximation of $x$ by elements of $C$, which is the closest point from $x$ to $C$. For uniformly convex and uniformly smooth Banach space $X$, $P_C: X \to C$ is a well-defined single-valued mapping.

In particular, if $X$ is a Hilbert space, then $P_C$ has the well-known basic variational principle: for any $x \in X$ and $u \in C$,

$$u = P_C x \quad \Leftrightarrow \quad \langle x - u, u - z \rangle \geq 0, \text{ for all } z \in C. \tag{1.2}$$

This basic variational principle of $P_C$ can be considered as the fundamental theorem of the analysis related to the metric projection operator in Hilbert spaces. It has played a very important role in optimization theory, variational inequality, fixed point problems in Hilbert spaces.

Another important property of the metric projection $P_C$ in Hilbert spaces is that $P_C$ is nonexpansive:
$$\|P_C x - P_C y\| \leq \|x - y\|, \text{ for any } x, y \in X. \tag{1.3}$$

One of the most important applications of the basic variational principle and the nonexpansive property of $P_C$ is to help define and study the Gâteaux type directional differentiability of $P_C$ (see [5–6, 9–11, 13, 16–17]. Some applications of the directional differentiability of $P_C$ in Hilbert spaces have been given (see [13–14]). In these papers, the basic variational principle and the nonexpansiveness of $P_C$ together provided strong tools for the studying of the directional differentiability of $P_C$.

In general, the metric projection $P_C$ in a uniformly convex and uniformly smooth Banach space $X$ does not enjoy the nonexpansive property (1.3). Let $J: X \to X^*$ be the normalized duality mapping in a uniformly convex and uniformly smooth Banach space $X$. $J$ is a well-defined single-valued mapping. With the help of $J$, the basic variational principle of $P_C$ becomes the following version: for any $x \in X$ and $u \in C$,
$$u = P_C x \iff \langle J(x-u), u-z \rangle \geq 0, \text{ for all } z \in C. \tag{1.4}$$

The normalized duality mapping $J$ in uniformly convex and uniformly smooth Banach spaces has many properties. However, $J$ is not a linear operator and has some non-convex properties (see [11]). It substantially increases the difficulties in studying the directional differentiability of $P_C$ in uniformly convex and uniformly smooth Banach spaces. In [8, 13–15, 7, 18–21, 23], some types of differentiability are studied for some operators in Banach spaces.

In this paper, we will define the Gâteaux directional differentiability of $P_C$ in uniformly convex and uniformly smooth Banach spaces (see Definition 4.1): for an arbitrary given $x \in X$ and for a vector $v \in X$ with $v \neq \theta$, the directional derivative of $P_C$ at point $x$ along direction $v$ is defined by

$$P'_C(x; v) = \lim_{t \downarrow 0} \frac{P_C(x+tv) - P_C(x)}{t}. \tag{1.5}$$

In section 4, we study the directional differentiability of $P_C$ onto nonempty closed and convex subsets in $X$. In sections 5 and 6, we investigate some properties of directional differentiability of $P_C$ onto closed balls and closed cones in $X$. In these two cases, we will provide the analytic representations (exact solutions) of $P_C$. For any $y \in C$, the inverse image of $y$ by the metric projection $P_C$ in $X$ is defined by

$$P_C^{-1}(y) = \{x \in X: P_C(x) = y\}.$$

In sections 5, 6, we will see that the inverse mapping $P_C^{-1}$ plays a crucial role in studying the directional differentiability of $P_C$ in uniformly convex and uniformly smooth Banach spaces.

Therefore, in section 3, we investigate some properties of $P_C^{-1}$, which will be used in the following sections 5 and 6.

In section 7, we introduce the concept of $p$-$q$ uniformly convex and uniformly smooth Banach spaces, and we will prove one of the main theorems in this paper (Theorem 7.2): If $X$ is a $p$-$q$ uniformly convex and uniformly smooth Banach space, then $P_C$ is directionally differentiable on the whole space $X$.

In section 2, we introduce the concepts of function and $J$-function of smoothness in smooth Banach spaces; these functions are denoted by $\psi$ and $\xi$, respectively. The $J$-function $\xi$ of smoothness is defined with respect to the normalized duality mapping $J$ in smooth Banach spaces. We establish a very close connection between $\psi$ and $\xi$, which implies that the smoothness of a Banach space can be described by the normalized duality mapping $J$. In section 5, we find that the representations of the directional derivatives of the metric projection onto closed balls in uniformly convex and uniformly smooth Banach spaces are closely related with the function of smoothness $\psi$ of the underlying Banach spaces.

Notice that by the definition (1.5), we consider the directional differentiability of $P_C$ at every point in uniformly convex and uniformly smooth Banach space $X$. Some authors only consider the directional differentiability of $P_C$ at the boundary of $C$ when $X$ is a Hilbert space (see [16, 21] and others).

## 2. Preliminaries

### 2.1. The function of smoothness of smooth Banach spaces

Since the theme of this paper is about the directional differentiability of the standard metric projection in Banach spaces, which describes the smoothness of the metric projection. It follows immediately that the smoothness of the considered Banach spaces will play an important role in this study. In this subsection, we recall the concepts and properties of smooth Banach spaces and introduce the definition of function of smoothness of smooth Banach spaces, which will be multiply used in this paper.

Let $(X, \|\cdot\|)$ be a real Banach space with topological dual space $(X^*, \|\cdot\|_*)$ (In this paper, all considered Banach spaces are real). Let $\langle \cdot, \cdot \rangle$ denote the real canonical evaluation pairing between $X^*$ and $X$. Let $S(X) = \{x \in X : \|x\| = 1\}$, which is the unit sphere in $X$. If, for every $x, v \in S(X)$, the following limit exists,

$$\lim_{t \downarrow 0} \frac{\|x+tv\| - \|x\|}{t},$$

then, $X$ is said to be smooth. This raises the following definition.

**Definition 2.1.** Let $X$ be a smooth Banach space. Define a function $\psi: S(X) \times S(X) \to \mathbb{R}_+$ by

$$\psi(x, v) = \lim_{t \downarrow 0} \frac{\|x+tv\| - \|x\|}{t}, \text{ for any } (x, v) \in S(X) \times S(X). \tag{1.1}$$

$\psi$ is called the function of smoothness of this smooth Banach space $X$.

If the limit in (1.1) is attained to $\psi(x, v)$ uniformly for $(x, v) \in S(X) \times S(X)$, then $X$ is said to be uniformly smooth (see Takahashi [22]). The modules of smoothness of $X$ is denoted by $\rho$, which is defined by

$$\rho(t) = \sup\left\{\frac{\|x+y\|+\|x-y\|}{2} - 1 : x, y \in X, \|x\| = 1, \|y\| = t\right\}, \text{ for } t > 0.$$

It $X$ is uniformly smooth, then the modules of smoothness $\rho$ of $X$ has the following properties:

(U$_1$) $\rho(t) \leq t$, for any $t > 0$;
(U$_2$) $\lim_{t \downarrow 0} \frac{\rho(t)}{t} = 0$;

### 2.2. The *J*-function of smoothness of smooth Banach spaces

The normalized duality mapping $J$ in Banach spaces play a very important role in the analysis of Banach spaces. Let $X$ be a Banach space with dual space $X^*$, in general, $J$ is a set-valued mapping $J: X \to 2^{X^*} \setminus \{\emptyset\}$, which is defined by

$$Jx = \{\varphi \in X^*: \langle \varphi, x \rangle = \|\varphi\|_* \|x\| = \|x\|^2 = \|\varphi\|_*^2\}, \text{ for any } x \in X.$$

We list the following properties of $J$ for easy reference (see Takahashi [22] for more details).

**Lemma 2.2.** *Let $X$ be a Banach space. Then*

(J$_1$) *$X$ is reflexive $\Rightarrow J: X \to X^*$ is onto;*
(J$_2$) *$X$ is smooth $\Leftrightarrow J$ is single-valued (Theorem 4.3.1 and Theorem 4.3.2);*
(J$_3$) *$X$ is smooth $\Rightarrow J$ is norm to weak\* continuous (Theorem 4.3.2 and Theorem 4.3.3);*
(J$_4$) *$X$ is uniformly smooth $\Rightarrow J: X \to X^*$ is norm to norm continuous (Theorem 4.3.5).*

The following proposition presents the connections between the smoothness and the normalized duality mapping in Banach spaces.

**Lemma 2.2.** *Let $X$ be a smooth Banach space with function of smoothness $\psi$. Then, for any $(x, v) \in S(X) \times S(X)$, the limit*

$$\xi(x, v) \equiv \lim_{t \downarrow 0} \frac{\langle J(x+tv), x \rangle - \langle J(x), x \rangle}{t}, \tag{1.2}$$

*exists such that*

$$\psi(x, v) = \frac{1}{2}(\langle J(x), v \rangle + \xi(x, v)), \text{ for any } (x, v) \in S(X) \times S(X). \tag{1.3}$$

The function $\xi: S(X) \times S(X) \to \mathbb{R}_+$ is called the *J*-function of smoothness of $X$.

*Proof.* For any $(x, v) \in S(X) \times S(X)$, by (1.1), we have

$$\psi(x, v) = \lim_{t \downarrow 0} \frac{\|x+tv\| - \|x\|}{t}$$
$$= \lim_{t \downarrow 0} \frac{\|x+tv\|^2 - \|x\|^2}{t(\|x+tv\| + \|x\|)}$$
$$= \lim_{t \downarrow 0} \frac{\langle J(x+tv),\ x+tv \rangle - \langle J(x),\ x \rangle}{t(\|x+tv\| + \|x\|)}$$
$$= \frac{1}{2\|x\|} \lim_{t \downarrow 0} \left( \frac{\langle J(x+tv),\ tv \rangle}{t} + \frac{\langle J(x+tv),\ x \rangle - \langle J(x),\ x \rangle}{t} \right)$$
$$= \frac{1}{2} \left( \lim_{t \downarrow 0} \langle J(x+tv),\ v \rangle + \lim_{t \downarrow 0} \frac{\langle J(x+tv),\ x \rangle - \langle J(x),\ x \rangle}{t} \right)$$

Since $x + tv \to x$, as $t \to 0$, in $X$ with respect to norm, by property (J3), we have

$$\lim_{t \downarrow 0} \langle J(x + tv),\ v \rangle = \langle J(x),\ v \rangle. \qquad \square$$

We list some notations below used in this paper. For any $u, v \in X$ with $u \neq v$, we write

(a) $\overline{v, u} = \{tv + (1-t)u \colon 0 \leq t \leq 1\}$;
(b) $\overrightarrow{v, u} = \{v + t(u - v) \colon 0 \leq t < \infty\}$;
(c) $\overleftrightarrow{v, u} = \{v + t(u - v) \colon -\infty < t < \infty\}$.

$\overline{v, u}$ is called a closed segment in $X$ with ending points at $u$ and $v$; $\overrightarrow{v, u}$ is called a closed ray in $X$ with ending point at $v$ and with direction $u - v$, which is a closed and convex cone with vertex at $v$, as a special case of cones in $X$; $\overleftrightarrow{v, u}$ is called a line in $X$ passing through points $v$ and $u$.

**Corollary 2.8 in [12].** *Let $X$ be a uniformly convex and uniformly smooth Banach space and $K$ a closed cone in $X$ with vertex at $v \in X$. Then*

(i) *If $v = \theta$, we have*
  (a) *$JK$ is a closed cone in $X^*$ with vertex at $\theta^* = J\theta$;*
  (b) *However, in general*
    $K$ *is convex* $\nRightarrow$ $JK$ *is convex;*
(ii) *If $v \neq \theta$ or $K$ is a ray with $\theta \notin K$, then $JK$ is not a cone (not a ray).*

**2.2. Metric projection in uniformly convex and uniformly smooth Banach Spaces**

Let $C$ be a nonempty closed and convex subset of the uniformly convex and uniformly smooth Banach space $X$. Let $P_C \colon X \to C$ denote the (standard) metric projection (operator). For any $x \in X$, $P_C x \in C$ such that

$$\|x - P_C x\| \leq \|x - z\|, \text{ for all } z \in C.$$

The metric projection operator has the following useful properties.

**Proposition 2.6 in [2].** *Let $X$ be a uniformly convex and uniformly smooth Banach space and $C$ a nonempty closed and convex subset of $X$. Then the metric projection $P_C \colon X \to C$ holds the following properties.*

(i)     *The operator $P_C$ is fixed on $C$; that is, $P_C(x) = x$, for any $x \in C$;*
(ii)    *$P_C$ has the following basic variational properties. For any given $x \in X$, for any $x \in X$ and $u \in C$, we have*

$$u = P_C(x) \iff \langle J_X(x - u), u - z \rangle \geq 0, \text{ for all } z \in C \text{ (see 5.e in [2])};$$
$$u = P_C(x) \iff \langle J_X(x - u), x - z \rangle \geq \|x - u\|^2, \text{ for all } z \in C \text{ (Theorem 5.1 in [2])};$$

(iii)   $\langle J_X(x - P_C(x)) - J_X(y - P_C(y)), P_C(x) - P_C(y) \rangle \geq 0$, *for any $x, y \in X$* (see 5.i in [2]);
(iv)   *In general, $P_C$ is not expansive* (see 5.f in [2]);
(v)    *$P_C$ is uniformly continuous on each bounded subset in $X$* (see 5.f in [2]);

## 3. Properties of inverse images of metric projection in uniformly convex and uniformly smooth Banach spaces

Throughout this section, let $(X, \|\cdot\|)$ be a real uniformly convex and uniformly smooth Banach space with topological dual space $(X^*, \|\cdot\|_*)$. Let $C$ be a nonempty closed and convex subset of $X$. In this section, we study the properties of the inverse images of points in $X$ by the metric projection $P_C$. These properties will be used in studying the differentiability of $P_C$ in the following sections.

**Definition 3.1.** Let $C$ be a nonempty closed and convex subset of $X$. For any $y \in C$, the inverse image of $y$ by the metric projection $P_C$ in $X$ is defined by

$$P_C^{-1}(y) = \{u \in X : P_C(u) = y\}.$$

**Lemma 3.2.** *Let $C$ be a nonempty closed and convex subset of $X$. Then,*

$$X = \bigcup_{y \in C} P_C^{-1}(y).$$

*Proof.* Since $X$ is a uniformly convex and uniformly smooth Banach space and $C$ is a nonempty closed and convex subset of $X$, it is well known that the metric projection $P_C : X \to C$ is a well-defined single-valued mapping from $X$ onto $C$. This lemma is proved.   □

**Theorem 3.1 in [12].** *Let $C$ be a nonempty closed and convex subset of $X$. For any $y \in C$, suppose that there is $x \in X \setminus C$ with $y = P_C(x)$. Then,*

(a) *$P_C^{-1}(y)$ is a closed cone with vertex at $y$ in $X$;*
(b) *However, in general, $P_C^{-1}(y)$ is not convex.*

**Definition 3.3.** Let $C$ be a nonempty closed and convex subset in $X$. Let $y \in C$.

(i)    If $P_C^{-1}(y) = \{y\}$, then $y$ is called an *internal point* of $C$;
(ii)   If $P_C^{-1}(y) \supsetneq \{y\}$, then $y$ is called a *cuticle point* of $C$.

The collection of all internal points of $C$ is denoted by $\mathbb{I}(C)$ and the collection of all cuticle points of $C$ is denoted by $\mathbb{C}(C)$. Then $\{\mathbb{I}(C), \mathbb{C}(C)\}$ forms a partition of $C$. More precisely, we have

$$\mathbb{I}(C) \cap \mathbb{C}(C) = \emptyset \quad \text{and} \quad \mathbb{I}(C) \cup \mathbb{C}(C) = C.$$

By Theorem 3.1 in [12], $\mathbb{I}(C)$ and $\mathbb{C}(C)$ have the following characteristics.

**Lemma 3.4**. *Let $C$ be a nonempty closed and convex subset of $X$. For any $y \in C$, we have*

(a) $y \in \mathbb{I}(C)$ *if and only if, for $u \in X$,*

$$y = P_C(u + y) \implies u = \theta;$$

(b) $y \in \mathbb{C}(C)$ *if and only if, there is $u \in X$ with $u \neq \theta$ such that*

$$y = P_C(u + y).$$

Next, we give some propositions to demonstrate the concepts of $\mathbb{I}(C)$ and $\mathbb{C}(C)$. For an arbitrary fixed $r > 0$, we respectively write the closed, open balls and sphere in $X$ with radius $r$ at center $c \in X$ by $B(c, r)$, $B^\circ(c, r)$ and $S(c, r)$. Then $B(c, r)$ is a closed and convex subset of $X$. In particular, when $c = \theta$ and $r = 1$, the unit closed ball and unit sphere are written by $B(\theta, 1) = B(X)$ and $S(\theta, 1) = S(X)$.

**Lemma 3.5**. *Let $c \in X$ and $r > 0$. Then we have*

$$P_{B(c,r)}(x) = c + \frac{r}{\|x-c\|}(x - c), \quad \text{for any } x \in X \backslash B(c, r). \tag{3.1}$$

*In particular, we have*

$$P_{B(X)}(x) = \frac{x}{\|x\|}, \quad \text{for any } x \in X \text{ with } \|x\| > 1.$$

*Proof.* For an arbitrary given $x \in X \backslash B(c, r)$, which is equivalent to $\|x - c\| > r$. Then, for any $z \in B(c, r)$, which satisfies $\|z - c\| \leq r$, by $1 - \frac{r}{\|x-c\|} > 0$, we have

$$\langle J(x - (c + \frac{r}{\|x-c\|}(x - c))), (c + \frac{r}{\|x-c\|}(x - c)) - z \rangle$$
$$= \left(1 - \frac{r}{\|x-c\|}\right) \langle J(x - c), \frac{r}{\|x-c\|}(x - c) + c - z \rangle$$
$$= \left(1 - \frac{r}{\|x-c\|}\right) (\langle J(x - c), \frac{r}{\|x-c\|}(x - c) \rangle + \langle J(x - c), c - z \rangle)$$
$$= \left(1 - \frac{r}{\|x-c\|}\right) (r\|x - c\| + \langle J(x - c), c - z \rangle)$$
$$\geq \left(1 - \frac{r}{\|x-c\|}\right) (r\|x - c\| - \|x - c\|\|c - z\|)$$
$$\geq 0, \qquad \text{for any } z \in B(c, r).$$

By the basic variational principle of $P_{B(c,r)}$, this proves (3.1). □

**Proposition 3.6.** *For any $c \in X$ and for an arbitrary fixed $r > 0$, we have*

(i)      $\mathbb{I}(B(c, r)) = B^o(c, r)$;
(ii)     $\mathbb{C}(B(c, r)) = S(c, r)$;
(iii)    *For any $y \in S(c, r)$,*

$$P^{-1}_{B(c,r)}(y) = \{y + t(y - c): 0 \leq t < \infty\}. \tag{3.1}$$

*Proof.* Take an arbitrary $y \in B^o(c, r)$. For any given $x \in X$ with $x \neq y$, there is $t > 0$, such that $y + t(x - y) \in B^o(c, r)$. Then

$$\langle J(x - y), y - (y + t(x - y))\rangle = \langle J(x - y), -t(x - y)\rangle = -t\|x - y\|^2 < 0.$$

By the basic variational principle of $P_{B(c,r)}$, the above inequality implies

$$x \notin P^{-1}_{B(c,r)}(y), \text{ for any } x \in X \text{ with } x \neq y.$$

Next, we prove (iii). It is clear to see that (iii) and (i) imply (ii).

Notice that $B(c, r) = B(\theta, r) + c$. For any $y \in S(c, r)$ any $t > 0$, we calculate

$$\begin{aligned}
&\langle J((y + t(y - c)) - y), y - z\rangle \\
&= t(\langle J(y - c), y - c + c - z\rangle) \\
&= t(\|y - c\|^2 + \langle J(y - c), c - z\rangle) \\
&\geq t(\|y - c\|^2 - \|y - c\|\|c - z\|) \\
&= t(r^2 - r\|c - z\|) \\
&\geq 0, \text{ for any } z \in B(c, r) \text{ (it satisfies } \|c - z\| \leq r\text{)}.
\end{aligned}$$

By the basic variational principle of $P_{B(c,r)}$ again, this implies

$$y + t(y - c) \in P^{-1}_{B(c,r)}(y), \text{ for any } t > 0.$$

This proves

$$\{y + t(y - c): 0 \leq t < \infty\} \subseteq P^{-1}_{B(c,r)}(y). \tag{3.2}$$

On the other hand, for any $x \in X$ with $\|x - c\| > r$, we take

$$w = c + \frac{r}{\|x-c\|}(x - c).$$

This implies that $w \in S(c, r)$ such that

$$x = c + \frac{\|x-c\|}{r}(w - c) = w + \left(\frac{\|x-c\|}{r} - 1\right)(w - c).$$

By $\frac{\|x-c\|}{r} - 1 > 0$, from the above proof, this implies

$$x \in P_{B(c,r)}^{-1}(w), \text{ for } w \in S(c, r).$$

From the proof of (3.2), this implies

$$x \in \{w + t(w - c): 0 \leq t < \infty\}, \text{ for } w \in S(c, r).$$

Since the metric projection $P_{B(c,r)}: X \to B(c, r)$ is one to one, for any $x \in X \setminus B(c, r)$, it implies

$$x \notin P_{B(c,r)}^{-1}(y) \implies x \notin \{y + t(y - c): 0 \leq t < \infty\}. \tag{3.3}$$

(3.1) follows from (3.2) and (3.3). □

**Proposition 3.7.** *Let $C$ be a proper closed subspace of $X$. Then, we have*

(i) $\mathbb{I}(C) = \emptyset$;
(ii) $\mathbb{C}(C) = C$.

*Proof.* Since $C$ is a proper closed subspace of $X$, by a separation theorem (a corollary of the Hahn-Banach Theorem, for example, see corollary 1.2.4 in Takahashi [22]), there is $\varphi \in X^*$ with $\|\varphi\|_* = 1$ such that

$$\langle \varphi, z \rangle = 0, \text{ for all } z \in C. \tag{3.4}$$

This implies $J^*\varphi \in X \setminus C$. For every given $y \in C$, since $J^*\varphi \in X \setminus C$ and $C$ is a proper closed subspace of $X$, it implies $J^*\varphi + y \notin C$. By (3.4) and the condition that $C$ is a closed subspace of $X$, we calculate

$$\langle J((J^*\varphi + y) - y), y - z \rangle$$
$$= \langle J(J^*\varphi), y - z \rangle$$
$$= \langle \varphi, y - z \rangle$$
$$= 0, \text{ for all } z \in C.$$

By the basic variational principle of $P_{B(v,r)}$, this implies

$$P_C(J^*\varphi + y) = y \text{ and } J^*\varphi + y \notin C, \text{ for every given } y \in C.$$

It proves this proposition. □

**Definition 3.8.** Let $C$ be a proper closed subspace of $X$. We denote

$$C^\perp = \{x \in X: \langle Jx, z \rangle = 0, \text{ for all } z \in C\}.$$

$C^\perp$ is a closed cone in $X$ with vertex at $\theta$ (by Corollary 2.8 in [12], $C^\perp$ is non-convex, in general). $C^\perp$ is called the orthogonal cone of $C$.

**Proposition 3.9.** *Let $C$ be a proper closed subspace of $X$. Then, we have*

(i) $P_C^{-1}(\theta) = C^\perp$;

(ii)    $y \in C \implies P_C^{-1}(y) = y + P_C^{-1}(\theta) = y + C^\perp.$

*Proof.* Proof of (i). Since $C$ is a proper closed subspace of $X$, by Proposition 3.7 and the basic variational principle of $P_C$, for $x \in X$, we have

$$x \in P_C^{-1}(\theta)$$
$$\Leftrightarrow \langle J(x - \theta), \theta - z \rangle \geq 0, \text{ for all } z \in C$$
$$\Leftrightarrow \langle Jx, z \rangle = 0, \text{ for all } z \in C$$
$$\Leftrightarrow x \in C^\perp.$$

Proof of (ii). For every $y \in C$, similarly to the proof of (i), for $x \in X$, we have

$$x \in P_C^{-1}(y)$$
$$\Leftrightarrow \langle J(x - y), y - z \rangle \geq 0, \text{ for all } z \in C$$
$$\Leftrightarrow \langle J(x - y), z \rangle = 0, \text{ for all } z \in C$$
$$\Leftrightarrow x - y \in C^\perp$$
$$\Leftrightarrow x \in y + C^\perp.$$
□

**Definition 3.10**. Let $K$ be a closed and convex cone in $X$ with vertex at $v \in X$. We denote

$$K^\wedge = P_K^{-1}(v) = \{x \in X : \langle J(x - v), v - z \rangle \geq 0, \text{ for all } z \in K\}.$$

$K^\wedge = P_K^{-1}(v)$ is a closed cone in $X$ with vertex at $v$ (by Theorem 3.1 in [12], $P_K^{-1}(v)$ is non-convex, in general). $K^\wedge$ is called the dual cone of $K$ in $X$.

**Lemma 3.11**. *Let $K$ be a closed and convex cone in $X$ with vertex at $v \in X$. Then, for every $y \in K$ with $y \neq v$ and for any $u = v + t(y - v)$ with $t > 0$, we have*

$$P_K^{-1}(u) = u - y + P_K^{-1}(y).$$

*Proof.* For $y \in K$, let $u = v + t(y - v)$, for $t > 0$. It follows that $u - y = (1 - t)(v - y)$. Then, for $x \in X$, we have

$$x + (u - y) \in P_K^{-1}(u)$$
$$\Leftrightarrow \langle J((x + (u - y)) - u), u - z \rangle \geq 0, \text{ for all } z \in K$$
$$\Leftrightarrow \langle J(x - y), v + t(y - v) - z \rangle \geq 0, \text{ for all } z \in K$$
$$\Leftrightarrow \langle J(x - y), (t - 1)(y - v) + y - z \rangle \geq 0, \text{ for all } z \in K$$
$$\Leftrightarrow \langle J(x - y), (t - 1)(y - v) \rangle = 0,$$
$$\quad \text{and } \langle J(x - y), y - z \rangle \geq 0, \text{ for all } z \in K$$
$$\Leftrightarrow x \in P_K^{-1}(y).$$

In the above proof, we applied the following implication

$$\langle J(x - y), (t - 1)(y - v) + y - z \rangle \geq 0, \text{ for all } z \in K.$$

$$\implies \langle J(x - y), (t - 1)(y - v) \rangle = 0 \text{ and } \langle J(x - y), y - z \rangle \geq 0, \text{ for all } z \in K.$$

This implication can be proved as follows. Take $z = y$, we get $\langle J(x - y), (t - 1)(y - v)\rangle \geq 0$; take $z = y + 2(t - 1)(y - v)$, we get $\langle J(x - y), -(t - 1)(y - v)\rangle \geq 0$. This proves

$$\langle J(x - y), (t - 1)(y - v)\rangle = 0. \qquad \square$$

**Corollary 3.12**. *Let $K$ be a closed and convex cone in $X$ with vertex at $\theta$. Then, for every $y \in K$ with $y \neq \theta$, we have*

$$y + P_K^{-1}(ty) = ty + P_K^{-1}(y), \text{ for any } t > 0.$$

*Proof.* This corollary follows by taking $u = ty$ in Lemma 3.11. $\qquad \square$

### 4. Directional differentiability of the metric projection onto closed and convex subsets

Similar to the previous section, in this section, let $X$ be a uniformly convex and uniformly smooth Banach space with dual space $X^*$ and $C$ a nonempty closed and convex subset of $X$. The metric projection $P_C$ is a well-defined single valued mapping from $X$ onto $C$, which assures that the definition of directionally differentiability of $P_C$ is meaningful. We define the Gâteaux directional differentiability of $P_C$ on the whole space $X$.

**Definition 4.1.** For $x \in X$ and $v \in X$ with $v \neq \theta$, if the following limit exists (that is a point in $X$),

$$\lim_{t \downarrow 0} \frac{P_C(x+tv) - P_C(x)}{t}$$

then, $P_C$ is said to be (Gâteaux) directionally differentiable at point $x$ along direction $v$, which is denoted by

$$P_C'(x; v) = \lim_{t \downarrow 0} \frac{P_C(x+tv) - P_C(x)}{t}.$$

$P_C'(x; v)$ is called the (Gâteaux) directional derivative of $P_C$ at point $x$ along direction $v$; and $v$ is called a (Gâteaux) differentiable direction of $P_C$ at $x$. If $P_C$ is (Gâteaux) directionally differentiable at point $x \in X$ along every direction $v \in X$ with $v \neq \theta$, then $P_C$ is said to be (Gâteaux) directionally differentiable at point $x \in X$. It is denoted by

$$P_C'(x)(v) = \lim_{t \downarrow 0} \frac{P_C(x+tv) - P_C(x)}{t}, \text{ for } v \in X \text{ with } v \neq \theta.$$

$P_C'(x)(v)$ is called the (Gâteaux) directional derivative of $P_C$ at point $x$ along direction $v$. Let $A$ be an open subset in $X$. If $P_C$ is (Gâteaux) directionally differentiable at every point $x \in A$, then $P_C$ is said to be (Gâteaux) directionally differentiable on $A \subseteq X$.

Notice that $P_C$ is (Gâteaux) directionally differentiable at point $x$ along direction $v$, if and only if, there is a point, denoted by $P_C'(x; v) \in X$, such that

$$P_C(x + tv) = P_C(x) + tP_C'(x; v) + o(t), \text{ for } t > 0.$$

By Definition 4.1, we give some properties of the directional differentiability of the metric projection $P_C$ onto nonempty closed and convex subsets of $X$.

**Proposition 4.2.** *Let $C$ be a nonempty closed and convex subset of $X$. Then, $P_C$ is directionally differentiable at every point $x \in X\backslash C$ along both directions $P_C(x) - x$ and $x - P_C(x)$, respectively, such that*
$$P'_C(x; x - P_C(x)) = P'_C(x; P_C(x) - x) = \theta, \text{ for every } x \in X\backslash C.$$

*Proof.* For any given point $x \in X\backslash C$, from Lemma 3.2, $P_C(x) \in C$ is well defined satisfying $x - P_C(x) \neq \theta$. By Theorem 3.1 in [12], $P_C^{-1}(P_C(x))$ is a closed cone with vertex at $P_C(x)$ in $X$, which yields
$$\overrightarrow{P_C(x), x} \subseteq P_C^{-1}(P_C(x)).$$

This implies
$$P_C(x) + t(x - P_C(x)) \in P_C^{-1}(P_C(x)), \text{ for all } t \geq 0.$$

That is,
$$P_C[P_C(x) + t(x - P_C(x))] = P_C(x), \text{ for all } t \geq 0.$$

We calculate
$$\begin{aligned}
& P'_C(x; x - P_C(x)) \\
&= \lim_{t \downarrow 0} \frac{P_C(x + t(x - P_C(x))) - P_C(x)}{t} \\
&= \lim_{t \downarrow 0} \frac{P_C(P_C(x) + (1+t)(x - P_C(x))) - P_C(x)}{t} \\
&= \lim_{t \downarrow 0} \frac{P_C(x) - P_C(x)}{t} \\
&= \theta.
\end{aligned}$$

We similarly calculate
$$\begin{aligned}
& P'_C(x; P_C(x) - x) \\
&= \lim_{t \downarrow 0} \frac{P_C(x + t(P_C(x) - x)) - P_C(x)}{t} \\
&= \lim_{t \downarrow 0} \frac{P_C(P_C(x) + (1-t)(x - P_C(x))) - P_C(x)}{t} \\
&= \lim_{t \downarrow 0, t \leq 1} \frac{P_C(P_C(x) + (1-t)(x - P_C(x))) - P_C(x)}{t} \\
&= \lim_{t \downarrow 0, t \leq 1} \frac{P_C(x) - P_C(x)}{t} \\
&= \theta.
\end{aligned}$$ □

**Lemma 4.3.** *Let $C$ be a nonempty closed and convex subset of $X$. Let $u, w \in C$ with $u \neq w$. For every point $x \in \overline{u, w}$, we have*

(a) *If $x \neq u$ and $x \neq w$, then $P_C$ is directionally differentiable at point $x$ along both directions $w - u$ and $u - w$ such that*
$$P'_C(x; u - w) = u - w \quad \text{and} \quad P'_C(x; w - u) = w - u;$$

(b) *$P_C$ is directionally differentiable at $u$ along both directions $w - u$ and $u - w$ such that*
$$P'_C(u; w - u) = w - u \quad \text{and} \quad P'_C(w; u - w) = u - w.$$

*Proof.* (a) Take any $x \in \overline{u, w}$ with $x \neq u$ and $x \neq w$. There is a positive number $\alpha$ with $0 < \alpha < 1$ such that $x = \alpha u + (1 - \alpha)w$. Since $C$ is convex, it yields $x \in \overline{u, w} \subseteq C$. We calculate

$$\begin{aligned}
P'_C(x; u - w) &= \lim_{t \downarrow 0} \frac{P_C(x + t(u-w)) - P_C(x)}{t} \\
&= \lim_{t \downarrow 0} \frac{P_C(w + (\alpha + t)(u-w)) - P_C(x)}{t} \\
&= \lim_{t \downarrow 0, t < 1-\alpha} \frac{P_C(w + (\alpha + t)(u-w)) - P_C(x)}{t} \\
&= \lim_{t \downarrow 0, t < 1-\alpha} \frac{w + (\alpha + t)(u-w) - x}{t} \\
&= \lim_{t \downarrow 0, t < 1-\alpha} \frac{t(u-w)}{t} \\
&= u - w.
\end{aligned}$$

In this proof, we used the fact when $0 < t < 1 - \alpha$, $w + (\alpha + t)(u - w) \in \overline{u, w} \subseteq C$. Rest of the proofs are similar to the above proof. □

**Lemma 4.4.** *Let $C$ be a nonempty closed and convex subset of $X$. Then, the following statements are equivalent*

(i)     $P_C$ *is directionally differentiable on $X$ such that, for every point $x \in X$,*

$$P'_C(x)(v) = \theta, \text{ for any } v \in X \text{ with } v \neq \theta;$$

(ii)    $P_C$ *is a constant operator; that is, $C$ is a singleton.*

*Proof.* (i) $\Rightarrow$ (ii). Assume, by the way of contradiction, that $C$ is not a singleton. Then, there are $u, w \in C$ with $u \neq w$. By parts (b) in Lemma 4.3, we have

$$P'_C(u)(w - u) = w - u \quad \text{and} \quad P'_C(w)(u - w) = u - w.$$

This contradicts to the assumption in part (i).

(ii) $\Rightarrow$ (i). Suppose that $C = \{y\}$, for some $y \in X$, which is a singleton. Then, $P_C$ is directionally differentiable on $X$. More precisely, for every point $x \in X$, we have

$$P'_C(x)(v) = \lim_{t \downarrow 0} \frac{P_C(x + t(u-w)) - P_C(x)}{t} = \lim_{t \downarrow 0} \frac{y - y}{t} = \theta, \text{ for any } v \in X \text{ with } v \neq \theta. \quad \square$$

**Corollary 4.5.** *Let $C$ be a nonempty closed and convex subset of $X$. Then, for every point $x \in X$, there is least one differentiable direction of $P_C$ at $x$.*

*Proof.* If $C$ is a singleton, then, this corollary follows from Lemma 4.4 immediately. Next, we suppose that $C$ is a nonempty, non-singleton, closed and convex subset of $X$. In this case, for every point $x \in X \setminus C$, from Lemma 4.2, $x$ has at least two differentiable directions $P_C(x) - x$ and $x - P_C(x)$; for any $x \in C$, since $C$ is not a singleton, there is $u \in C$ with $u \neq x$. By part (c) of Lemma 4.3, $x$ has at least one differentiable direction $u - x$. □

Next lemma proves that the directional derivative of $P_C$ is positive homogenous.

**Lemma 4.6.** *Let $C$ be a nonempty closed and convex subset of $X$. For $x \in X$ and for $v \in X$ with $v \neq \theta$, if $P_C$ is directionally differentiable at $x$ along direction $v$, then, for any $\lambda > 0$, $P_C$ is directionally differentiable at $x$ along direction $\lambda v$ such that*

$$P_C'(x; \lambda v) = \lambda P_C'(x; v), \text{ for any } \lambda > 0.$$

*Proof.* For $x \in X$ and for $v \in X$ with $v \neq \theta$, suppose that $P_C'(x; v) = \lim_{t \downarrow 0} \frac{P_C(x+tv) - P_C(x)}{t}$ exists. Then, for any $\lambda > 0$, we have

$$\begin{aligned} P_C'(x; \lambda v) &= \lim_{t \downarrow 0} \frac{P_C(x+t\lambda v) - P_C(x)}{t} \\ &= \lambda \lim_{t \downarrow 0} \frac{P_C(x+t\lambda v) - P_C(x)}{\lambda t} \\ &= \lambda P_C'(x; v). \end{aligned}$$
□

**Proposition 4.7.** *Let $C$ be a nonempty closed and convex subset of $X$. Let $y \in C$. Suppose $(P_C^{-1}(y))^o \neq \emptyset$. Then, $P_C$ is directionally differentiable on $(P_C^{-1}(y))^o$ such that, for any $x \in (P_C^{-1}(y))^o$, we have*

$$P_C'(x)(v) = \theta, \text{ for every } v \in X \text{ with } v \neq \theta.$$

*Proof.* For any given point $x \in (P_C^{-1}(y))^o$, there is $\delta_x > 0$ such that

$$x + w \in (P_C^{-1}(y))^o, \text{ for any } w \in X \text{ with } \|w\| < \delta_x.$$

That is,

$$P_C(x+w) = P_C(x) = y, \text{ for any } w \in X \text{ with } \|w\| < \delta_x. \quad (4.1)$$

For any $v \in X$ with $v \neq \theta$, by (4.1), we calculate

$$\begin{aligned} P_C'(x)(v) &= \lim_{t \downarrow 0} \frac{P_C(x+tv)) - P_C(x)}{t} \\ &= \lim_{t \downarrow 0, t < \delta_x \backslash \|v\|} \frac{P_C(x+tv)) - y}{t} \\ &= \lim_{t \downarrow 0, t < \delta_x \backslash \|v\|} \frac{P_C(x) - y}{t} \\ &= \lim_{t \downarrow 0, t < \delta_x \backslash \|v\|} \frac{y - y}{t} \\ &= \theta. \end{aligned}$$
□

**Proposition 4.8.** *Let $C$ be a nonempty closed and convex subset of $X$. Suppose $C^o \neq \emptyset$. Then $P_C$ is directionally differentiable on $C^o$ such that, for any $x \in C^o$, we have*

$$P_C'(x)(v) = v, \quad \text{for every } v \in X \text{ with } v \neq \theta.$$

*Proof.* The proof of this proposition is similar to the proof of Proposition 4.7. For any given point $x \in C^o$, there is $\delta_x > 0$ such that

$$x + w \in C^o, \text{ for any } w \in X \text{ with } \|w\| < \delta_x.$$

This implies

$$P_C(x) = x \text{ and } P_C(x+w) = x+w, \text{ for any } w \in X \text{ with } \|w\| < \delta_x. \quad (4.2)$$

For any $v \in X$ with $v \neq \theta$, by (4.2), we calculate

$$\begin{aligned} P_C'(x)(v) &= \lim_{t \downarrow 0} \frac{P_C(x+tv)) - P_C(x)}{t} \\ &= \lim_{t \downarrow 0, t < \delta_x \setminus \|v\|} \frac{P_C(x+tv)) - x}{t} \\ &= \lim_{t \downarrow 0, t < \delta_x \setminus \|v\|} \frac{x+tv - x}{t} \\ &= \lim_{t \downarrow 0, t < \delta_x \setminus \|v\|} \frac{tv}{t} \\ &= v. \end{aligned}$$

□

## 5. Directional differentiability of the metric projection onto closed balls

As what we defined in section 4, in this section, let $X$ be a uniformly convex and uniformly smooth Banach space with dual space $X^*$. Before we study the directional differentiability of the metric projection on closed balls in $X$, we need to introduce the following notations. Recall that, for $c \in X$ and $r > 0$, the closed, open balls and sphere in $X$ with center at $c$ and with radius $r$ are respectively written as $B(c, r)$, $(B(c, r))^o$ and $S(c, r)$.

For any $x \in S(c, r)$, we define two subsets $x_{(c,r)}^{\uparrow}$ and $x_{(c,r)}^{\downarrow}$ of $X \setminus \{\theta\}$ as follows: for $v \in X$ with $v \neq \theta$, we say

(a) $v \in x_{(c,r)}^{\uparrow} \iff$ there is $\delta > 0$ such that $\|(x+tv) - c\| \geq r$, for all $t \in (0, \delta)$;
(b) $v \in x_{(c,r)}^{\downarrow} \iff$ there is $\delta > 0$ such that $\|(x+tv) - c\| < r$, for all $t \in (0, \delta)$.

In particular, $B(X)$ is the closed unit ball and $S(X)$ is the closed unit sphere. For any given $x \in S(X)$ and for $v \in X$ with $v \neq \theta$, we write

(c) $x^{\uparrow} = x_{(\theta,1)}^{\uparrow}$: $v \in x^{\uparrow} \iff$ there is $\delta > 0$ such that $\|x + tv\| \geq 1$, for all $t \in (0, \delta)$;
(d) $x^{\downarrow} = x_{(\theta,1)}^{\downarrow}$: $v \in x^{\downarrow} \iff$ there is $\delta > 0$ such that $\|x + tv\| < 1$, for all $t \in (0, \delta)$.

Next lemma shows that, for any given $x \in S(c, r)$, the two subsets $x_{(c,r)}^{\uparrow}$ and $x_{(c,r)}^{\downarrow}$ form a partition of $X \setminus \{\theta\}$.

**Lemma 5.1.** *Let $c \in X$ and $r > 0$. Then, for any $x \in S(c, r)$, we have*

$$x_{(c,r)}^{\uparrow} \cap x_{(c,r)}^{\downarrow} = \emptyset \quad \text{and} \quad x_{(c,r)}^{\uparrow} \cup x_{(c,r)}^{\downarrow} = X \setminus \{\theta\}.$$

*Proof.* For any $x \in S(X)$, by definitions, the disjoint property of $x_{(c,r)}^{\uparrow}$ and $x_{(c,r)}^{\downarrow}$ is evident. We only prove the second part. For any $v \in X$ with $v \neq \theta$, if $v \notin x_{(c,r)}^{\uparrow}$, then, there is $\delta > 0$, such that

$$\|(x + \delta v) - c\| < r. \tag{5.1}$$

It follows that, for any $t$ with $0 < t < \delta$, by (5.1), we have

$$\begin{aligned}&\|(x + tv) - c\| \\ &= \left\|\left(1 - \tfrac{t}{\delta}\right)(x - c) + \tfrac{t}{\delta}(x + \delta v - c)\right\| \\ &\leq \left(1 - \tfrac{t}{\delta}\right)\|x - c\| + \tfrac{t}{\delta}\|(x + \delta v) - c\| \\ &< r, \text{ for any } t \in (0, \delta).\end{aligned}$$

This implies $v \in x^{\downarrow}_{(c,r)}$. □

In the following theorem, we prove that the metric projection on closed balls is differentiable on the whole space $X$. Furthermore, we give the exact analytic representations (solutions) of the derivatives of the metric projection on closed balls.

**Theorem 5.2.** *Let $C = B(c, r)$ be a closed ball in $X$. Then, $P_C$ is directionally differentiable on $X$ such that, for every $v \in X$ with $v \neq \theta$, we have*

(i) *For any $x \in (B(c, r))^o$,*

$$P'_C(x)(v) = v;$$

(ii) *For any $x \in X \backslash B(c, r)$,*

$$P'_C(x)(v) = \frac{r}{\|x-c\|^2}\left(\|x - c\|v - \psi(\tfrac{x-c}{\|x-c\|}, \tfrac{v}{\|v\|})\|v\|(x - c)\right);$$

(iii) *For any $x \in S(c, r)$, we have*

(a) $$P'_C(x)(v) = v - \tfrac{\|v\|}{r}\psi(\tfrac{x-c}{\|x-c\|}, \tfrac{v}{\|v\|})(x - c), \text{ if } v \in x^{\uparrow}_{(c,r)};$$

(b) $$P'_C(x)(v) = v, \text{ if } v \in x^{\downarrow}_{(c,r)}.$$

*Proof.* It is clear that the directional differentiability of $P_C$ on the whole space $X$ follows immediately from the parts (i–iii). Part (i) is a consequence of Proposition 4.8. We prove part (ii). For any $x \in X \backslash B(c, r)$, and $v \in X$ with $v \neq \theta$, there is $\delta_x > 0$ such that

$$x + w \in X \backslash B(c, r), \text{ for any } w \in X \text{ with } \|w\| < \delta_x.$$

That is,

$$\|x - c\| > r \text{ and } \|(x + w) - c\| > r, \text{ for any } w \in X \text{ with } \|w\| < \delta_x.$$

By (3.1) in Lemma 3.5, we have

$$P_C(x) = c + \frac{r(x-c)}{\|x-c\|}$$

and

$$P_C(x + w) = c + \frac{r(x+w-c)}{\|x+w-c\|}, \text{ for any } w \in X \text{ with } \|w\| < \delta_x.$$

Since $X$ is a uniformly convex and uniformly smooth Banach space, then, the function of smoothness of $X$, $\psi: S(X)\times S(X) \to \mathbb{R}_+$ satisfies that

$$\lim_{t\downarrow 0}\frac{\|y+tw\| - \|y\|}{t} = \psi(y,w), \text{ uniformly for any } (y,w) \in S(X)\times S(X). \tag{5.2}$$

This implies that

$$P'_C(x)(v) = \lim_{t\downarrow 0}\frac{P_C(x+tv)-P_C(x)}{t}$$

$$= \lim_{t\downarrow 0, t<\delta_x\backslash\|v\|} \frac{\left(c+\frac{r(x+tv-c)}{\|x+tv-c\|}\right)-\left(c+\frac{r(x-c)}{\|x-c\|}\right)}{t}$$

$$= r\lim_{t\downarrow 0, t<\delta_x\backslash\|v\|} \frac{\frac{x+tv-c}{\|x+tv-c\|}-\frac{x-c}{\|x-c\|}}{t}$$

$$= r\lim_{t\downarrow 0, t<\delta_x\backslash\|v\|} \frac{\|x-c\|(x+tv-c) - \|x+tv-c\|(x-c)}{t\|x-c\|\|x+tv-c\|}$$

$$= r\lim_{t\downarrow 0, t<\delta_x\backslash\|v\|} \frac{\|x-c\|tv - (\|x+tv-c\|-\|x-c\|)(x-c)}{t\|x-c\|\|x+tv-c\|}$$

$$= r\lim_{t\downarrow 0, t<\delta_x\backslash\|v\|} \frac{v}{\|x+tv-c\|} - r\frac{x-c}{\|x-c\|^2}\lim_{t\downarrow 0, t<\delta_x\backslash\|v\|}\frac{\|x-c+tv\|-\|x-c\|}{t}$$

$$= \frac{rv}{\|x-c\|} - \frac{r(x-c)}{\|x-c\|^2}\lim_{t\downarrow 0, t<\delta_x\backslash\|v\|}\frac{\|x-c\|\left(\left\|\frac{x-c}{\|x-c\|}+t\frac{\|v\|}{\|x-c\|}\frac{v}{\|v\|}\right\|-\left\|\frac{x}{\|x-c\|}\right\|\right)}{t}$$

$$= \frac{rv}{\|x-c\|} - \frac{r(x-c)}{\|x-c\|^2}\lim_{t\downarrow 0, t<\delta_x\backslash\|v\|}\frac{\|v\|}{\|x-c\|}\cdot\frac{\|x-c\|\left(\left\|\frac{x-c}{\|x-c\|}+t\frac{\|v\|}{\|x-c\|}\frac{v}{\|v\|}\right\|-\left\|\frac{x-c}{\|x-c\|}\right\|\right)}{t\frac{\|v\|}{\|x-c\|}}$$

$$= \frac{rv}{\|x-c\|} - \frac{r\|v\|(x-c)}{\|x-c\|^2}\lim_{t\downarrow 0, t<\delta_x\backslash\|v\|}\frac{\left\|\frac{x-c}{\|x-c\|}+t\frac{\|v\|}{\|x-c\|}\frac{v}{\|v\|}\right\|-\left\|\frac{x-c}{\|x-c\|}\right\|}{t\frac{\|v\|}{\|x-c\|}}$$

$$= \frac{rv}{\|x-c\|} - \frac{r\|v\|(x-c)}{\|x-c\|^2}\psi\left(\frac{x-c}{\|x-c\|},\frac{v}{\|v\|}\right)$$

$$= \frac{r}{\|x-c\|^2}\left(\|x-c\|v - \psi\left(\frac{x-c}{\|x-c\|},\frac{v}{\|v\|}\right)\|v\|(x-c)\right).$$

Next, we prove (a) of (iii). Take any $x \in S(c, r)$ and any $v \in X$ with $v \neq \theta$ and $v \in x^\uparrow_{(c,r)}$. It follows that there is $\delta > 0$ such that $\|x + tv - c\| \geq r$, for all $t \in (0, \delta)$. By (3.1) and (5.2), we calculate

$$P'_C(x)(v) = \lim_{t\downarrow 0}\frac{P_C(x+tv)-P_C(x)}{t}$$

$$= \lim_{t\downarrow 0, t<\delta}\frac{\left(c+\frac{r}{\|x+tv-c\|}(x+tv-c)\right)-\left(c+\frac{r}{\|x-c\|}(x-c)\right)}{t}$$

$$= r \lim_{t\downarrow 0, t<\delta} \frac{\frac{x+tv-c}{\|x+tv-c\|} - \frac{x-c}{\|x-c\|}}{t}$$

$$= r \lim_{t\downarrow 0, t<\delta} \frac{\frac{tv}{\|x+tv-c\|}}{t} + r(x-c) \lim_{t\downarrow 0, t<\delta} \frac{\frac{1}{\|x+tv-c\|} - \frac{1}{\|x-c\|}}{t}$$

$$= \frac{rv}{\|x-c\|} - \frac{r(x-c)}{\|x-c\|^2} \lim_{t\downarrow 0, t<\delta} \frac{\|x-c+tv\| - \|x-c\|}{t}$$

$$= \frac{rv}{\|x-c\|} - \frac{r(x-c)}{\|x-c\|^2} \lim_{t\downarrow 0, t<\delta} \frac{\|x-c\|\left(\left\|\frac{x-c}{\|x-c\|} + t\frac{\|v\|}{\|x-c\|}\frac{v}{\|v\|}\right\| - \left\|\frac{x-c}{\|x-c\|}\right\|\right)}{t}$$

$$= \frac{rv}{\|x-c\|} - \frac{r(x-c)}{\|x-c\|^2} \lim_{t\downarrow 0, t<\delta} \frac{\|x-c\|\frac{\|v\|}{\|x-c\|}\left(\left\|\frac{x-c}{\|x-c\|} + t\frac{\|v\|}{\|x-c\|}\frac{v}{\|v\|}\right\| - \left\|\frac{x-c}{\|x-c\|}\right\|\right)}{t\frac{\|v\|}{\|x-c\|}}$$

$$= \frac{rv}{\|x-c\|} - \frac{r\|v\|(x-c)}{\|x-c\|^2} \lim_{t\downarrow 0, t<\delta} \frac{\left\|\frac{x-c}{\|x-c\|} + t\frac{\|v\|}{\|x-c\|}\frac{v}{\|v\|}\right\| - \left\|\frac{x-c}{\|x-c\|}\right\|}{t\frac{\|v\|}{\|x-c\|}}$$

$$= \frac{rv}{\|x-c\|} - \frac{r\|v\|(x-c)}{\|x-c\|^2} \psi\left(\frac{x-c}{\|x-c\|}, \frac{v}{\|v\|}\right)$$

$$= v - \frac{\|v\|}{r} \psi\left(\frac{x-c}{\|x-c\|}, \frac{v}{\|v\|}\right)(x-c).$$

Finally, we prove (b) of (iii). Take any $x \in S(c, r)$ and any $v \in X$ with $v \neq \theta$. Suppose $v \in x^{\downarrow}_{(c,r)}$. Then, there is $\delta > 0$ such that $\|x + tv - c\| < r$, for all $t \in (0, \delta)$. That is, $x \in B(c, r)$ and $x + tv \in (B(c, r))^o$, for all $t \in (0, \delta)$. Then, by Proposition 3.5, we have

$$P'_C(x)(v) = \lim_{t\downarrow 0} \frac{P_C(x+tv)) - P_C(x)}{t}$$
$$= \lim_{t\downarrow 0, t<\delta} \frac{x+tv - x}{t}$$
$$= v. \qquad \square$$

**Corollary 5.3.** *Let X be a uniformly convex and uniformly smooth Banach space. Then, $P_{B(X)}$ is directionally differentiable on X such that, for every $v \in X$ with $v \neq \theta$, we have*

(i) *For any $x \in (B(X))^o$,*

$$P'_{B(X)}(x)(v) = v;$$

(ii) *For any $x \in X \setminus B(X)$,*

$$P'_{B(X)}(x)(v) = \frac{1}{\|x\|^2}\left(\|x\|v - \psi\left(\frac{x}{\|x\|}, \frac{v}{\|v\|}\right)\|v\|x\right);$$

(iii) *For any $x \in S(X)$, we have*

(a) $\quad P'_{B(X)}(x)(v) = v - \psi(x, \frac{v}{\|v\|}) \|v\| x, \text{ if } v \in x^{\uparrow};$

(b) $\quad P'_{B(X)}(x)(v) = v, \text{ if } v \in x^{\downarrow}.$

## 6. Directional differentiability of the metric projection onto subspaces and cones

Throughout this section, let $X$ be a uniformly convex and uniformly smooth Banach space with dual space $X^*$. In this section, we study the directional differentiability of the metric projection on subspaces and cones of $X$. We also give the analytic representations of the directional derivatives of the metric projection.

**Theorem 6.1**. *Let $C$ be a proper closed subspace of $X$. Then, $P_C$ is directionally differentiable at every point $y \in C$ along any direction $v \in C^\perp \setminus \{\theta\}$ such that*

$$P'_C(y; v) = \theta, \text{ for any } v \in C^\perp \text{ with } v \neq \theta.$$

*Proof.* For the given proper closed subspace $C$ of $X$ in this theorem, by Definition 3.9, the orthogonal cone of $C$ is

$$C^\perp = \{x \in X: \langle Jx, z \rangle = 0, \text{ for all } z \in C\}.$$

$C^\perp$ is a closed cone in $X$ with vertex at $\theta$. For any $y \in C$, by Proposition 3.9, we have

$$P_C^{-1}(y) = y + P_C^{-1}(\theta) = y + C^\perp.$$

$P_C^{-1}(y)$ is also a closed cone in $X$ with vertex at $y$, which is non-convex, in general (by Corollary 2.8 in [12]). It follows that, for any given $y \in C$ and for any $v \in C^\perp$ with $v \neq \theta$, we have

$$y + tv \in P_C^{-1}(y), \text{ for ant } t > 0.$$

This implies

$$\begin{aligned} P'_C(y; v) &= \lim_{t \downarrow 0} \frac{P_C(y+tv)) - P_C(y)}{t} \\ &= \lim_{t \downarrow 0} \frac{P_C(y) - P_C(y)}{t} \\ &= \theta. \end{aligned}$$

$\square$

Before we study the directional differentiability of the metric projection on cones in uniformly convex and uniformly smooth Banach spaces, we recall the following example, which has been used in [11].

**Example 6.2**. Let $X = \mathbb{R}^3$ equipped with the 3-norm $\|\cdot\|_3$ defined, for any $z = (z_1, z_2, z_3) \in \mathbb{R}^3$, by

$$\|z\|_3 = \sqrt[3]{|z_1|^3 + |z_2|^3 + |z_3|^3}.$$

Then, $(\mathbb{R}^3, \|\cdot\|_3)$ is a uniformly convex and uniformly smooth Banach space (and it is not a Hilbert space). The dual space of $(\mathbb{R}^3, \|\cdot\|_3)$ is $(\mathbb{R}^3, \|\cdot\|_{\frac{3}{2}})$, which has the same underlying space $\mathbb{R}^3$ with the $\|\cdot\|_{\frac{3}{2}}$-norm such that, for $\varphi = (\varphi_1, \varphi_2, \varphi_3) \in \mathbb{R}^3$,

$$\|\varphi\|_{\frac{3}{2}} = \left(|\varphi_1|^{\frac{3}{2}} + |\varphi_2|^{\frac{3}{2}} + |\varphi_3|^{\frac{3}{2}}\right)^{\frac{2}{3}}.$$

The normalized duality mappings $J$ and $J^*$ on $(\mathbb{R}^3, \|\cdot\|_3)$ and $(\mathbb{R}^3, \|\cdot\|_{\frac{3}{2}})$, respectively, satisfy the following conditions: For any $z = (z_1, z_2, z_3) \in (\mathbb{R}^3, \|\cdot\|_3)$ with $z \neq \theta$,

$$Jz = \left(\frac{|z_1|^2 \text{sign}(z_1)}{\|z\|_3}, \frac{|z_2|^2 \text{sign}(z_2)}{\|z\|_3}, \frac{|z_3|^2 \text{sign}(z_3)}{\|z\|_3}\right) \in (\mathbb{R}^3, \|\cdot\|_{\frac{3}{2}}). \tag{6.1}$$

For any $\varphi = (\varphi_1, \varphi_2, \varphi_3) \in (\mathbb{R}^3, \|\cdot\|_{\frac{3}{2}})$ with $\varphi \neq \theta^*$,

$$J^*\varphi = \left(\frac{|\varphi_1|^{\frac{3}{2}-1} \text{sign}(\varphi_1)}{\left(\|\varphi\|_{\frac{3}{2}}\right)^{\frac{3}{2}-2}}, \frac{|\varphi_2|^{\frac{3}{2}-1} \text{sign}(\varphi_2)}{\left(\|\varphi\|_{\frac{3}{2}}\right)^{\frac{3}{2}-2}}, \frac{|\varphi_3|^{\frac{3}{2}-1} \text{sign}(\varphi_3)}{\left(\|\varphi\|_{\frac{3}{2}}\right)^{\frac{3}{2}-2}}\right) \in (\mathbb{R}^3, \|\cdot\|_3). \tag{6.2}$$

In next example, we consider a special closed and convex cone $K$ in $(\mathbb{R}^3, \|\cdot\|_3)$, which is called the positive cone of $(\mathbb{R}^3, \|\cdot\|_3)$. In order to investigate the directional differentiability of $P_K$, in the following example, we study the analytic representation of $P_K$.

**Proposition 6.3**. Let $K$ be the positive cone of uniformly convex and uniformly smooth Banach space $(\mathbb{R}^3, \|\cdot\|_3)$:

$$K = \{z = (z_1, z_2, z_3) \in \mathbb{R}^3 : z_i \geq 0, \text{ for } i = 1, 2, 3\}.$$

$K$ is a closed and convex cone in $(\mathbb{R}^3, \|\cdot\|_3)$. Let $K^o$ and $\partial K$ denote the interior and boundary of $K$, respectively. Take distinct numbers $i, j, k$ with $\{i, j, k\} = \{1, 2, 3\}$. Then, for any $z = (z_1, z_2, z_3) \in \mathbb{R}^3$, it is formally rewritten as $[z_i, z_j, z_k]$.

(i) For any $[y_i, y_j, y_k] \in K^o$, that is, $y_i > 0$, $y_j > 0$ and $y_k > 0$, we have

$$P_K^{-1}([y_i, y_j, y_k]) = [y_i, y_j, y_k];$$

(ii) For $[y_i, y_j, y_k] \in \partial K$ with $y_i > 0$, $y_j > 0$, and $y_k = 0$, we have

$$P_K^{-1}([y_i, y_j, 0]) = \{[y_i, y_j, x_k] \in \mathbb{R}^3 : x_k \leq 0\},$$

which is a ray in $\mathbb{R}^3$ with ending point at $[y_i, y_j, 0]$ and passing through $[y_i, y_j, -1]$;

(iii) For $[y_i, y_j, y_k] \in \partial K$ with $y_i > 0$, $y_j = 0$ and $y_k = 0$, we have

$$P_K^{-1}([y_i, 0, 0]) = \{[y_i, x_j, x_k] \in \mathbb{R}^3 : x_j \leq 0, x_k \leq 0\},$$

which is a closed and convex cone in $\mathbb{R}^3$ with vertex at $[z_i, 0, 0]$, which has its boundary enclosed by the following two rays:

(a) with ending point at $[y_i, 0, 0]$ and passing through $[y_i, 0, -1]$,
(b) with ending point at $[y_i, 0, 0]$ and passing through $[y_i, -1, 0]$;

(iv)
$$P_K^{-1}(\theta) = K^{\perp} = \{x = (x_1, x_2, x_3) \in \mathbb{R}^3 : x_i \leq 0, \text{ for } i = 1, 2, 3\}.$$

*Proof.* Parts (i) and (iv) are evident. The proofs of (ii) and (iii) are similar. So, we only prove (ii). For any $[y_i, y_j, x_k] \in \mathbb{R}^3$ with $x_k < 0$, by (6.1), we have

$$\langle J([y_i, y_j, x_k] - [y_i, y_j, 0]), [y_i, y_j, 0] - [z_i, z_j, z_k] \rangle$$
$$= \langle J([0, 0, x_k]), [y_i, y_j, 0] - [z_i, z_j, z_k] \rangle$$
$$= \langle [0, 0, x_k], [y_i, y_j, 0] - [z_i, z_j, z_k] \rangle$$
$$= -x_k z_k$$
$$\geq 0, \text{ for all } [z_i, z_j, z_k] \in K.$$

By the basic variational principle of $P_K$, this implies

$$[y_i, y_j, x_k] \in P_K^{-1}([y_i, y_j, 0]), \text{ for any } x_k < 0. \tag{6.3}$$

On the other hand, for any $[x_i, x_j, x_k] \in \mathbb{R}^3$ with $[x_i, x_j] \neq [y_i, y_j]$, then, one can check that there is $[z_i, z_j, z_k] \in K$ such that

$$\langle J([x_i, x_j, x_k] - [y_i, y_j, 0]), [y_i, y_j, 0] - [z_i, z_j, z_k] \rangle < 0.$$

By the basic variational principle of $P_K$ again, this implies

$$[x_i, x_j, x_k] \notin P_K^{-1}([y_i, y_j, 0]), \text{ if } [x_i, x_j] \neq [y_i, y_j]. \tag{6.4}$$

By (6.3) and (6.4), part (ii) is proved. □

In next proposition, we investigate the analytic representations of the directional derivatives of the metric projection $P_K$, in which $K$ is the positive cone defined in Proposition 6.3. Through the following proposition, we see that the analytic representations of $P_K$ is complicated, even though the positive $K$ is considered to be simple.

**Proposition 6.4.** Let $K$ be the positive cone of $(\mathbb{R}^3, \|\cdot\|_3)$. Then, $P_K$ is directionally differentiable on the whole space $(\mathbb{R}^3, \|\cdot\|_3)$ such that, any $v \in \mathbb{R}^3$ with $v \neq \theta$,

(I)   $y = [y_i, y_j, y_k] \in K$.

(i)   For every $y \in K^o$, we have
$$P'_K(y)(v) = v;$$

(ii) For every $y = [y_i, y_j, 0] \in \partial K$ with $y_i > 0$, $y_j > 0$, we have

$$P'_K(y)(v) = \begin{cases} v, & \text{for } v_k \geq 0 \\ [v_i, v_j, 0], & \text{for } v_k < 0. \end{cases}$$

(iii) For every $y = [y_i, 0, 0] \in \partial K$ with $y_i > 0$ we have

$$P'_K(y)(v) = \begin{cases} v, & \text{for } v_j \geq 0, v_k \geq 0, \\ [v_i, v_j, 0], & \text{for } v_j \geq 0, v_k < 0, \\ [v_i, 0, 0], & \text{for } v_j < 0, v_k < 0. \end{cases}$$

(iv) $$P'_K(\theta)(v) = \begin{cases} v, & \text{for } v_i \geq 0, v_j \geq 0, v_k \geq 0, \\ [v_i, v_j, 0], & \text{for } v_i \geq 0, v_j \geq 0, v_k < 0, \\ [v_i, 0, 0], & \text{for } v_i \geq 0, v_j < 0, v_k < 0, \\ [0, 0, 0], & \text{for } v_i < 0, v_j < 0, v_k < 0. \end{cases}$$

(II) $x = [x_i, x_j, y_k] \in \mathbb{R}^3 \setminus K$.

(v) For every $x = [x_i, x_j, y_k]$ with $x_i > 0$, $x_j > 0$, $x_k < 0$, we have

$$P'_K(x)(v) = [v_i, v_j, 0];$$

(vi) For every $x = [x_i, x_j, y_k]$ with $x_i > 0$, $x_j = 0$, $x_k < 0$, we have

$$P'_K(x)(v) = \begin{cases} [v_i, v_j, 0], & \text{for } v_j \geq 0, v_k < 0, \\ [v_i, 0, 0], & \text{for } v_j < 0, v_k < 0. \end{cases}$$

(vii) For every $x = [x_i, x_j, y_k]$ with $x_i = 0$, $x_j = 0$, $x_k < 0$, we have

$$P'_K(y)(v) = \begin{cases} [v_i, v_j, 0], & \text{for } v_i \geq 0, v_j \geq 0, \\ [v_i, 0, 0], & \text{for } v_i \geq 0, v_j < 0, \\ [0, 0, 0], & \text{for } v_i < 0, v_j < 0. \end{cases}$$

(viii) For every $x = [x_i, x_j, y_k]$ with $x_i < 0$, $x_j = 0$, $x_k < 0$, we have

$$P'_K(x)(v) = \begin{cases} [0, v_j, 0], & \text{for } v_j \geq 0, \\ [0, 0, 0], & \text{for } v_j < 0. \end{cases}$$

(ix) For every $x = [x_i, x_j, y_k]$ with $x_i < 0$, $x_j < 0$, $x_k < 0$, we have

$$P'_K(x)(v) = \theta.$$

*Proof.* Proof of (i). It immediately follows from Proposition 4.8. Then we prove parts (ii) and (ii). With the similar proofing ideas, parts (iv−ix) can be similarly proved.

Proof of (ii). For $t > 0$, since $y + tv = [y_i + tv_i, y_j + tv_j, tv_k]$. There is $\delta > 0$ such that $y_i + tv_i > 0$ and $y_j + tv_j > 0$, for any $t \in (0, \delta]$. By part (ii) of Proposition 6.3, for any $t \in (0, \delta]$, we have

$$P_K(y + tv) = \begin{cases} [y_i + tv_i, y_j + tv_j, tv_k], & \text{for } v_k \geq 0, \\ [y_i + tv_i, y_j + tv_j, 0], & \text{for } v_k < 0. \end{cases}$$

This implies

$$P'_K(y)(v) = \lim_{t \downarrow 0, t < \delta} \frac{P_K(y+tv)) - P_K(y)}{t}$$

$$= \begin{cases} \lim_{t \downarrow 0, t < \delta} \frac{[y_i + tv_i, y_j + tv_j, tv_k] - y}{t}, & \text{for } v_k \geq 0 \\ \lim_{t \downarrow 0, t < \delta} \frac{[y_i + tv_i, y_j + tv_j, 0] - y}{t}, & \text{for } v_k < 0, \end{cases}$$

$$= \begin{cases} \lim_{t \downarrow 0, t < \delta} \frac{tv}{t}, & \text{for } v_k \geq 0 \\ \lim_{t \downarrow 0, t < \delta} \frac{[tv_i, tv_j, 0]}{t}, & \text{for } v_k < 0, \end{cases}$$

$$= \begin{cases} v, & \text{for } v_k \geq 0 \\ [v_i, v_j, 0], & \text{for } v_k < 0. \end{cases}$$

Proof of (iii). For $t > 0$, since $y + tv = [y_i + tv_i, tv_j, tv_k]$. There is $\delta > 0$ such that $y_i + tv_i > 0$, for any $t \in (0, \delta]$. By parts (ii) and (iii) of Proposition 6.3, for any $t \in (0, \delta]$, we have

$$P_K(y + tv) = \begin{cases} [y_i + tv_i, tv_j, tv_k], & \text{for } v_j \geq 0, v_k \geq 0, \\ [y_i + tv_i, tv_j, 0], & \text{for } v_j \geq 0, v_k < 0, \\ [y_i + tv_i, 0, 0], & \text{for } v_j < 0, v_k < 0 \end{cases}$$

This implies

$$P'_K(y)(v) = \lim_{t \downarrow 0, t < \delta} \frac{P_K(y+tv)) - P_K(y)}{t}$$

$$= \begin{cases} \lim_{t \downarrow 0, t < \delta} \frac{[y_i + tv_i, tv_j, tv_k] - y}{t}, & \text{for } v_j \geq 0, v_k \geq 0, \\ \lim_{t \downarrow 0, t < \delta} \frac{[y_i + tv_i, tv_j, 0] - y}{t}, & \text{for } v_j \geq 0, v_k < 0, \\ \lim_{t \downarrow 0, t < \delta} \frac{[y_i + tv_i, 0, 0] - y}{t}, & \text{for } v_j < 0, v_k < 0 \end{cases}$$

$$= \begin{cases} \lim_{t\downarrow 0, t<\delta} \dfrac{tv}{t}, & \text{for } v_j \geq 0, v_k \geq 0, \\ \lim_{t\downarrow 0, t<\delta} \dfrac{[tv_i,\ tv_j,\ 0]}{t}, & \text{for } v_j \geq 0, v_k < 0, \\ \lim_{t\downarrow 0, t<\delta} \dfrac{[tv_i,\ 0,\ 0]}{t}, & \text{for } v_j < 0, v_k < 0 \end{cases}$$

$$= \begin{cases} v, & \text{for } v_j \geq 0, v_k \geq 0, \\ [v_i,\ v_j,\ 0], & \text{for } v_j \geq 0, v_k < 0, \\ [v_i,\ 0,\ 0], & \text{for } v_j < 0, v_k < 0. \end{cases}$$

The proofs of parts (iv−viii) are similar to the proofs of parts (ii, iii) and they are omitted here. By part (iv) of Proposition 6.3, notice that

$$\{[x_i, x_j, y_k] \in \mathbb{R}^3 : x_i < 0, x_j < 0, x_k < 0\} = (P_K^{-1}(\theta))^o.$$

By Proposition 4.7, this implies part (ix) immediately. □

### 7. $p$-$q$ uniformly convex and uniformly smooth Banach spaces

Similar to the previous sections, throughout this section, let $X$ be a uniformly convex and uniformly smooth Banach space with dual space $X^*$. Let $\delta$ and $\rho$ denote the modulus of convexity and modulus of smoothness of $X$, respectively. Both $\delta$ and $\rho$ have been used to study the geometric structures of uniformly convex and uniformly smooth Banach spaces (see [3,17, 22, 23]). They have the following well-known properties:

(a) Both $\delta(\varepsilon)$ and $\dfrac{\delta(\varepsilon)}{\varepsilon}$ are non-decreasing, for $\varepsilon \in (0, 2]$;

(b) $\rho(t) \leq t$, for $t > 0$ and $\dfrac{\rho(t)}{t} \to 0$, as $t \downarrow 0$.

The <u>Pisier</u> renorming theorem states that a super-reflexive uniformly convex and uniformly smooth Banach space $X$ admits an equivalent uniformly convex and uniformly smooth norm for which the modulus of convexity and the modulus of smoothness $\delta$ and $\rho$ satisfy the following conditions: there are positive numbers $a$, $b$ with $a \geq 1, b \geq 1$ and $p, q$ with $1 < p < q$ such that

(c) $\delta(\varepsilon) \geq a\varepsilon^p$, for $\varepsilon \in (0, 2]$;
(d) $\rho(t) \leq bt^q$, for $t > 0$.

This positive number $p > 1$ describes the "level" of convexity of $X$ and this number $q > 1$ describes the "level" of smoothness of $X$. The above observation motivates us to introduce the following definition.

**Definition 7.1.** Let $X$ be a uniformly convex and uniformly smooth Banach space. If there are positive numbers $a, b$ with $a \geq 1, b \geq 1$ and $p, q$ with $1 < p < q$ such that $\delta$ and $\rho$ satisfy the above conditions (c) and (d), then $X$ is called a $p$-$q$ uniformly convex and uniformly smooth Banach space.

Let $C$ be a nonempty closed and convex subset of $X$. It is well known that, if $X$ is a Hilbert space, then, the metric projection $P_C$ is nonexpansive. From the properties of $P_C$ recalled in section 2, we have that, in general uniformly convex and uniformly smooth Banach spaces, $P_C$ is not nonexpansive (see 5.e *in* [2]). However, $P_C$ is uniformly continuous on each bounded subset in $X$ (see 5.f *in* [2]).

In [1–3, 7], the authors investigate the estimation of the distance function related to $P_C$ by using the modulus of convexity and the modulus of smoothness $\delta$ and $\rho$ of $X$.

**Theorems** in [1–3, 7]. *Let $X$ be a uniformly convex and uniformly smooth Banach space and $C$ a nonempty closed and convex subset of $X$. For any $x, y \in X$ with $x \neq y$, let*

$$k = 2\max\{1, \|x - P_C y\|, \|P_C x - y\|\}. \tag{7.1}$$

*Then, we have*

$$\|P_C x - P_C y\| \leq k\delta^{-1}(6\rho(2\|x - y\|)), \tag{5.4 in [1]}$$

$$\|P_C x - P_C y\| \leq k\delta^{-1}(\rho(8kL\|x - y\|)). \tag{5.5 in [1]}$$

Now, we state and prove the main theorem in this section.

**Theorems 7.2**. *Let $X$ be a $p$-$q$ uniformly convex and uniformly smooth Banach space and $C$ a nonempty closed and convex subset of $X$. Then $P_C$ is directionally differentiable on $X$.*

*Proof.* We need to show that, for any $x \in X$, the following limit exists:

$$\lim_{t \downarrow 0} \frac{P_K(x+tv)) - P_K(x)}{t}, \text{ for any } v \in X \text{ with } v \neq \theta. \tag{7.2}$$

For this given $p$-$q$ uniformly convex and uniformly smooth Banach space $X$, there are 4 positive numbers $a, b$ with $a \geq 1$, $b \geq 1$ and $p, q$ with $1 < p < q$ satisfying conditions (c) and (d).

Then, for any given fixed $v \in X$ with $v \neq \theta$, take a positive number $\lambda < 1$ such that

$$(6a^{-1}b(2\|v\|\lambda)^q)^{\frac{1}{p}} < 2. \tag{7.3}$$

For any numbers $s, t$ with $0 < s < t < \lambda < 1$, we have

$$\left\| \frac{P_C(x+tv)) - P_C(x)}{t} - \frac{P_C(x+sv)) - P_C(x)}{s} \right\|$$

$$\leq \frac{1}{t}\|P_C(x+tv) - P_C(x+sv)\| + \frac{t-s}{st}\|P_C(x+sv) - P_C(x)\|. \tag{7.4}$$

By the uniformly continuity of $P_C$ on each bounded subset in $X$, one can see that

$$\sup\{\|P_C(x+tv)-(x+sv)\|: \text{for } 0<s<t<1\}<\infty,$$
$$\sup\{\|(x+tv)-P_C(x+sv)\|: \text{for } 0<s<t<1\}<\infty,$$
$$\sup\{\|P_C(x+sv)-x\|: \text{for } 0<s<1\}<\infty,$$
$$\sup\{\|(x+sv)-P_Cx\|: \text{for } 0<s<1\}<\infty.$$

Hence, in (7.1), we take $k>0$ as follows

$$k = 2\max\{1, \sup\{\|P_C(x+tv)-(x+sv)\|: \text{for } 0<s<t<1\},$$
$$\sup\{\|(x+tv)-P_C(x+sv)\|: \text{for } 0<s<t<1\},$$
$$\sup\{\|P_C(x+sv)-x\|: \text{for } 0<s<1\},$$
$$\sup\{\|(x+sv)-P_Cx\|: \text{for } 0<s<1\}\}. \tag{7.5}$$

By substituting (7.5) into (5.4) in [1] listed in Theorems in [1−3, 7], and by (c), (d) and (7.4), for any $0<s<t<\lambda$, we have

$$\left\|\frac{P_C(x+tv))-P_C(x)}{t} - \frac{P_C(x+sv))-P_C(x)}{s}\right\|$$

$$\leq \frac{k}{t}\delta^{-1}(6\rho(2(t-s)\|v\|)) + k\frac{t-s}{st}\delta^{-1}(6\rho(2s\|v\|))$$

$$\leq \frac{k}{t}\delta^{-1}(6b(2(t-s)\|v\|)^q) + k\frac{t-s}{st}\delta^{-1}(6b(2s\|v\|)^q)$$

$$\leq \frac{k}{t}\delta^{-1}(a6a^{-1}b(2\|v\|t)^q) + \frac{k}{s}\delta^{-1}(a6a^{-1}b(2\|v\|s)^q)$$

$$= \frac{k}{t}\delta^{-1}\left(a\left((6a^{-1}b(2\|v\|t)^q)^{\frac{1}{p}}\right)^p\right) + \frac{k}{s}\delta^{-1}\left(a\left((6a^{-1}b(2\|v\|s)^q)^{\frac{1}{p}}\right)^p\right)$$

By (7.3), $0<\lambda<1$ satisfying

$$0 < (6a^{-1}b(2\|v\|s)^q)^{\frac{1}{p}} < (6a^{-1}b(2\|v\|t)^q)^{\frac{1}{p}} < (6a^{-1}b(2\|v\|\lambda)^q)^{\frac{1}{p}} < 2.$$

By (a), $\delta^{-1}$ is non-decreasing. We have

$$\frac{k}{t}\delta^{-1}\left(a\left((6a^{-1}b(2\|v\|t)^q)^{\frac{1}{p}}\right)^p\right) + \frac{k}{s}\delta^{-1}\left(a\left((6a^{-1}b(2\|v\|s)^q)^{\frac{1}{p}}\right)^p\right)$$

$$\leq \frac{k}{t}(6a^{-1}b(2\|v\|t)^q)^{\frac{1}{p}} + \frac{k}{s}(6a^{-1}b(2\|v\|s)^q)^{\frac{1}{p}}$$

$$= k(6a^{-1}b(2\|v\|)^q)^{\frac{1}{p}} t^{\frac{q}{p}-1} + k(6a^{-1}b(2\|v\|)^q)^{\frac{1}{p}} s^{\frac{q}{p}-1}. \tag{7.6}$$

From (7.6), we have

$$\left\|\frac{P_C(x+tv))-P_C(x)}{t} - \frac{P_C(x+sv))-P_C(x)}{s}\right\|$$

$$\leq k(6a^{-1}b(2\|v\|)^q)^{\frac{1}{p}} t^{\frac{q}{p}-1} + k(6a^{-1}b(2\|v\|)^q)^{\frac{1}{p}} s^{\frac{q}{p}-1}, \text{ for any } s, t \text{ with } 0 < s < t < \lambda.$$

Notice that $\frac{q}{p} - 1 > 0$ and for this arbitrarily given $v \in X$ with $v \neq 0$, $k(6a^{-1}b(2\|v\|)^q)^{\frac{1}{p}}$ is a constant, which is independent from $1 > \lambda > 0$ and $s, t$ with $0 < s < t < \lambda < 1$. This implies that

$$\left\|\frac{P_C(x+tv))-P_C(x)}{t} - \frac{P_C(x+sv))-P_C(x)}{s}\right\|$$

can be arbitrarily small as what one desires, for any $s, t$ with $0 < s < t < \lambda < 1$, as $\lambda$ is chosen to be small enough. This proves (7.2) and it completes the proof of this theorem. □

**Corollary 7.3**. *Let X be a p-q uniformly convex and uniformly smooth Banach space and C a nonempty closed and convex subset of X. Then, for any $x \in X$,*

$$\lim_{t\downarrow 0} \frac{P_K(x+tv))-P_K(x)}{t} = P'_K(x)(v), \text{ uniformly for } v \in S(X). \tag{7.2}$$

*Proof.* In the proof of Theorem 7.2, for any $v \in S(X)$, we have

$$\left\|\frac{P_C(x+tv))-P_C(x)}{t} - \frac{P_C(x+sv))-P_C(x)}{s}\right\|$$

$$\leq k(6a^{-1}b2^q)^{\frac{1}{p}} t^{\frac{q}{p}-1} + k(6a^{-1}b2^q)^{\frac{1}{p}} s^{\frac{q}{p}-1}, \text{ for any } s, t \text{ with } 0 < s < t < \lambda.$$

This implies this corollary. □

### 8. Applications to Hilbert spaces

Since every Hilbert space is a special uniformly convex and uniformly smooth Banach space, all results obtained in the previous sections are hold for Hilbert spaces, in which, representations for the directional derivatives in Hilbert spaces should be simpler. We start at the function of smoothness of Hilbert spaces.

**Lemma 8.1**. *Let X be a Hilbert space. Then, for all $(x, v) \in S(X) \times S(X)$, we have*

$$\psi(x, v) = \psi(v, x) = \xi(x, v) = \xi(v, x) = \langle x, v \rangle. \tag{8.1}$$

*Proof.* For any $(x, v) \in S(X) \times S(X)$, we calculate

$$\psi(x, v) = \lim_{t\downarrow 0} \frac{\|x+tv\| - \|x\|}{t}$$
$$= \lim_{t\downarrow 0} \frac{\|x+tv\|^2 - \|x\|^2}{t(\|x+tv\| + \|x\|)}$$

$$= \lim_{t\downarrow 0} \frac{\langle x+tv,\ x+tv\rangle - \langle x,\ x\rangle}{t(\|x+tv\| + \|x\|)}$$

$$= \frac{1}{2\|x\|} \lim_{t\downarrow 0} \frac{2t\langle x,\ v\rangle + t^2}{t}$$

$$= \langle x, v\rangle.$$

We can similarly prove that $\xi(x, v) = \langle x, v\rangle$. □

As a consequence of Theorem 8.2, we have

**Corollary 8.2.** *Let $C = B(c, r)$ be a closed ball in a Hilbert space X. Then, $P_C$ is directionally differentiable on X such that, for every $v \in X$ with $v \neq \theta$, we have*

(j) *For any $x \in (B(c, r))^o$,*

$$P'_C(x)(v) = v;$$

(ii) *For any $x \in X\setminus B(c, r)$,*

$$P'_C(x)(v) = \frac{r}{\|x-c\|^3}(\|x - c\|^2 v - \langle x, v\rangle(x - c));$$

(iii) *For any $x \in S(c, r)$, we have*

(a) $$P'_C(x)(v) = v - \frac{1}{r\|x-c\|}\langle x, v\rangle(x - c), \ if\ v \in x^{\uparrow}_{(c,r)};$$

(b) $$P'_C(x)(v) = v, \ if\ v \in x^{\downarrow}_{(c,r)}.$$

*Proof.* Proof of (ii). By part (ii) of Theorem 5.2 and Lemma 8.1, we have

$$P'_C(x)(v) = \frac{r}{\|x-c\|^2}\left(\|x - c\|v - \psi(\frac{x-c}{\|x-c\|}, \frac{v}{\|v\|})\|v\|(x - c)\right)$$

$$= \frac{r}{\|x-c\|^2}\left(\|x - c\|v - \frac{1}{\|x-c\|\|v\|}\langle x, v\rangle\|v\|(x - c)\right)$$

$$= \frac{r}{\|x-c\|^3}(\|x - c\|^2 v - \langle x, v\rangle(x - c)).$$

Proof of (b) in (iii).

$$P'_C(x)(v) = v - \frac{\|v\|}{r}\frac{1}{\|x-c\|\|v\|}\langle x, v\rangle(x - c)$$

$$= v - \frac{1}{r\|x-c\|}\langle x, v\rangle(x - c).$$ □

**Corollary 8.3.** *Let X be a Hilbert space. Then, $P_{B(X)}$ is directionally differentiable on X such that, for every $v \in X$ with $v \neq \theta$, we have*

(i) *For any $x \in (B(X))^o$,*
$$P'_{B(X)}(x)(v) = v;$$

(ii) *For any $x \in X \backslash B(X)$,*
$$P'_{B(X)}(x)(v) = \frac{1}{\|x\|^3}(\|x\|^2 v - \langle x, v \rangle x),$$
*and*
$$P'_{B(X)}(x)(x) = \theta;$$

(iii) *For any $x \in S(X)$,*
   (a) $P'_{B(X)}(x)(v) = v - \langle x, v \rangle x, \ if \ v \in x^\uparrow$,
   (b) $P'_{B(X)}(x)(v) = v, \ if \ v \in x^\downarrow$;
   (c) $x \in S(X) \Longrightarrow x \in x^\uparrow$,
$$P'_{B(X)}(x)(x) = \theta, \ if \ x \in S(X).$$

Proof. It follows from Corollary 8.2 immediately. □

## Acknowledgments

The author is very grateful to Professors Phil Blau, Akhtar Khan for their kind communications in the development stage of this paper.